\documentclass[12pt]{amsart}
\usepackage{amscd,amssymb,epsfig}
\calclayout

\newif\ifpdf

\ifx\pdfoutput\undefined

  \pdffalse 
\else 
  \pdftrue\pdfoutput=1 
\fi 
 
\numberwithin{equation}{section} 
\def\RR{\mathbb R}

\newcommand{\ud}{\,d} 
\newcommand{\Sym}{\mathbb{S}} 
\newcommand{\R}{\mathbb{R}}

\newcommand\M{\mathbb{M}}
\renewcommand\P{{\mathcal P}}

\newcommand\tr{\operatorname{tr}}

\newcommand\x{\times}
\renewcommand{\div}{\operatorname{div}} 
\newcommand{\grad}{\operatorname{grad}}

\newcommand{\curl}{\operatorname{curl}}

\newcommand{\as}{\operatorname{as}} 
\newcommand{\tir}[1]{\ensuremath{\overline {#1}}} 
\newtheorem{thm}{Theorem}[section]

\def\whsq{\vbox to 5.8pt 
{\offinterlineskip\hrule 
\hbox to 5.8pt{\vrule height 
5.1pt\hss\vrule height 5.1pt}\hrule}}

\def\<{\langle} 
\def\>{\rangle} 
\def\PP{{\mathop{{\rm I}\kern-.2em{\rm P}}\nolimits}} 
\def\FF{{\mathop{{\rm I}\kern-.2em{\rm F}}\nolimits}}   
\def\ZZ{{\mathop{{\rm I}\kern-.2em{\rm Z}}\nolimits}} 
 
\newlength{\sidemargin} 
\begin{document}
\title[]{
Rectangular Mixed Elements for Elasticity with Weakly Imposed symmetry Condition
}

\thanks{The author was supported in part by
NSF grant DMS-0811052 and the Sloan Foundation. Fruitful discussions with D. Arnold are gratefully acknowledged}

\author{Gerard Awanou}
\address{Department of Mathematical Sciences,
Northern Illinois University,
Dekalb, IL, 60115}
\email{awanou@math.niu.edu}
\urladdr{http://www.math.niu.edu/\char'176awanou}

\begin{abstract}
We present new rectangular mixed finite elements for linear elasticity. 
The approach is based on a modification of the
Hellinger-Reissner functional in which the symmetry of the stress field is enforced
weakly through the introduction of a Lagrange multiplier.
The elements are analogues of the lowest order elements described in
Arnold, Falk and Winther [
Mixed finite element methods for linear elasticity with weakly imposed symmetry.
Mathematics of Computation 76 (2007), pp. 1699--1723]. 
Piecewise constants are used to approximate the
displacement and the rotation. 
The first order BDM elements are used to approximate each row of the stress field.

\end{abstract}

\maketitle
\section{Introduction}
The theory of elasticity is used to predict the 
response of a material to applied forces.
The unknowns in the equations are the stress field, a symmetric matrix 
field which encodes the internal forces and 
the displacement, a vector field. For various reasons, mixed finite elements where one approximates both the 
stress and displacement are the methods of choice. 
One seeks the stress in the space of symmetric matrix fields 
with components  
square integrable and with divergence, taken row-wise, also square integrable. The displacement
is sought in the space of square integrable vector fields. The pair forms a unique saddle
point of the
Hellinger-Reissner functional. It is very difficult to construct at the discrete level, finite element 
spaces which satisfy Brezzi's stability conditions. These conditions 
provide sufficient conditions
for the stability of mixed finite element methods. Indeed for several decades before the work of Arnold and Winter \cite{Arnold2002,Arnold2003}
the existence of such elements was an open problem.
These elements have been extended to rectangular meshes in two dimension \cite{Arnold2005,Chen2010}, three dimension \cite{Awanou10} and
on tetrahedral meshes \cite{Arnold2008,Adams2005}.
Despite their relative complexity, mixed finite elements with symmetric stress fields
are useful in certain situations \cite{Nicaise2008}. If one desires simpler elements, one is forced to turn
to nonconforming elements. Nonconformity can be introduced by weakening the symmetry condition or by
weakening the requirement that the stress field is $L^2$ integrable.
We refer to \cite{Awanou10b} for a review on nonconforming elements with symmetric stress fields and other approaches to linear
elasticity.

Stable mixed finite elements with weakly imposed symmetry have been introduced in
\cite{Amara1979,Arnold1984a,Stenberg1986,Stenberg1988,Stenberg1988a,Morley1989,Arnold2006,Arnold2007,Falk08,Guzman10c,Guzman10b,Guzman10d},
The purpose of this paper is to present elements with weakly imposed symmetry for rectangular meshes.
Precisely, we will  use 
piecewise constants to approximate the
displacement and the rotation and 18 or 12 dimensional spaces to approximate the stress field.
The simplest older element on rectangular meshes in two dimensions is the one of \cite{Morley1989} with 11 degrees of freedom for the stress, piecewise
constants to approximate the displacement but a 4 dimensional space to approximate the rotation. The advantage of our element
is that the rotation can be eliminated by static condensation. In three dimensions as well, our elements are simpler than Morley's elements.

The paper is organized as follows: after some preliminaries in the next section, 
we present our low order elements in two dimension and then in three dimension. We conclude with
some remarks on higher order elements. 

\section{Preliminaries}
Let $\Omega$ be a simply connected polygonal domain of $\RR^n, n=2,3$, occupied by a
linearly elastic body which is clamped on $\partial \Omega$.
We denote as usual by $L^2 (\Omega, \RR^n)$ the space of square
integrable vector fields with values in $\RR^n$ and $H^k (K, X)$ the
space of functions with domain $K \subset \RR^n$, taking values in the finite
dimensional space $X$, and with all derivatives of order at most $k$
square integrable. 
We let $H(\div,
\Omega,X)$
be the space of square-integrable fields taking values in $X$ 
and which have square integrable
divergence. 
For our purposes, $X$ will be either $\mathbb{M}$ the space of 
$n \times n$ matrices, $\Sym$ the space of $n \times n$ symmetric matrices, $\RR^n,
\mbox{ or } \RR$, and in the latter case, we simply write $H^k(X)$. 
The divergence operator is the usual divergence for vector fields which produces a matrix field
when acting on a matrix field by taking the divergence of each row. 
We will also need the space $H(\curl,\Omega,\RR^n)$ of square-integrable fields with
square integrable $\curl$. We recall that in two dimension for a scalar function $q$,
 $\curl(q)  = (\partial_2 q, -\partial_1 q) $ and 
in three dimension 
\begin{align*}
\curl(q_1,q_2,q_3) & = (\partial_2 q_3-\partial_3 q_2, -\partial_1 q_3+\partial_3 q_1, \partial_1 q_2-
\partial_2 q_1).
\end{align*}
For a vector field in two dimension or a matrix field in three dimension, the curl operator
produces a matrix field by taking the curl of each row.
The
norms in $H^k(K,X)$ and $H^k(K)$ are denoted respectively by $|| \cdot
||_{H^k}$ and $|| \cdot ||_k$. We use the usual notations of $\mathcal{P}_k(K,X)$
for the space of polynomials on $K$ with values in $X$ of total degree
less than $k$ and $\mathcal{P}_{k_1, k_2}(K,X)$ for the space of polynomials of
degree at most $k_1$ in $x$ and of degree at most $k_2$ in $y$. 
Similarly, 
$\mathcal{P}_{k_1, k_2, k_3}(K,X)$ denotes the space of polynomials of
degree at most $k_1$ in $x$, of degree at most $k_2$ in $y$ and 
of degree at most $k_3$ in $z$.
We
write $ \mathcal{P}_k $, $\mathcal{P}_{k_1, k_2}$ 
and $\mathcal{P}_{k_1, k_2, k_3}$ respectively when $X = \RR$. 

The solution $(\sigma, u) \in H (\div, \Omega, \Sym) \times L^2 (\Omega, \RR^n)$
of the elasticity problem can be characterized as the unique critical point
of the Hellinger-Reissner functional 
$$
\mathcal{J} (\sigma, v) = \int_{\Omega} \left( \dfrac{1}{2} A \tau + \div
\tau \cdot v - f \cdot v \right) \, d x .
$$
The compliance tensor $A = A (x) : \Sym \to \Sym$ is given, bounded and symmetric
positive definite uniformly with respect to $x \in \Omega$, and the body
force $f$ is also given. In the homogeneous and isotropic case,
$$
A \sigma = \frac{1}{2 \mu} \bigg( \sigma - \frac{\lambda}{2 \mu + 2 \lambda} \mathrm{tr} \ (\sigma) I \bigg)
$$
where $I$ is the identity matrix and $\lambda$ and $\mu$ are the positive Lame constants.

To treat both two and three dimensional problems in a unified framework, one possibility is to use finite element differential forms
\cite{Arnold2007b}. However, for $n=2,3$ a simple device will suffice. We define
$\mathbb{P}$ to be $\R$ when $n=2$ and $\mathbb{P}=\R^3$ for $n=3$. Then we define 
$\as \tau=\tau_{12}-\tau_{21}$ for a $2 \times 2$ matrix and 
$\as{\tau}=(\tau_{32}-\tau_{23},\tau_{13}-\tau_{31},\tau_{21}-\tau_{12})'$ in three dimension.
For a symmetric matrix field, $\as \tau =0$. 
Next, we define
$\mathbb{H}$ to be $\R^2$ when $n=2$ and $\mathbb{H}=\M$ for $n=3$.
For the formulation with weakly imposed symmetry condition, a critical point
of the extended functional
$$
\mathcal{J} (\sigma, v) + \int_{\Omega} \eta \cdot \as \tau 
$$
is sought over $H(\div,\Omega,\mathbb{M}) 
\times L^2(\Omega,\mathbb{R}^n) \times L^2(\Omega,\mathbb{P})$. 
The unique solution $(\sigma, u, \gamma)$ satisfies
\begin{align}
\begin{split} \label{elasSys}
(A \sigma,\tau) + (\div \tau,u) + (\as \tau,\gamma)& = 0, \quad \tau \in  H(\div,\Omega,\mathbb{M}),\\
(\div \sigma,v) & =(f,v),  \quad v \in L^2(\Omega,\mathbb{R}^n),\\
(\as \sigma,q) & = 0, \quad q \in L^2(\Omega,\mathbb{P}).
\end{split}
\end{align}
 For the associated discrete system with finite element spaces $\Sigma_h\times V_h \times Q_h 
 \subset H(\div,\Omega,\mathbb{M}) \times L^2(\Omega,\mathbb{R}^n)\times L^2(\Omega,\mathbb{P})$, the symmetry
 condition will be enforced only weakly. The Brezzi's conditions for stability are
 \begin{itemize}
\item There exists a positive constant $c_1$ independent of $h$ such that
$||\tau||_{H({\div})} \leq c_1 (A \tau,\tau) $, if $\tau \in \Sigma_h$,
$  (\div \tau, v) =0$ for all $v \in V_h$ and $  ( \as \tau, q) =0, \forall q \in 
Q_h$,
\item There exists a positive constant $c_2$ independent of $h$ such that $\forall \ (v,q) \in V_h \times Q_h, (v,q) \neq (0,0), 
\exists \
\tau \in \Sigma_h, \tau \neq 0$ with $ (\div \tau, v)  + ( \as \tau, q)  \geq 
c_2  
||\tau||_{H(\div)} (||v||_{L^2} + ||q||_{L^2}) $.
\end{itemize}

To fulfill these conditions, we construct $\Sigma_h, V_h$ and $Q_h$ such that
\begin{enumerate}
\item[1-] $\div \Sigma_h \subset V_h$
\item[2-] Given $(v,q) \in V_h \times Q_h, (v,q) \neq (0,0), \exists \, \tau \in \Sigma_h, 
\tau \neq 0$ such that 
\begin{equation}
||\tau||_{H(\div)} \leq C (||v||_{L^2} + ||q||_{L^2}), \label{divrel}
\end{equation}
and $\div \tau =v$, $P_{Q_h} \as \tau=q$, where $P_{Q_h}$ is the $L^2$ projection operator. 
\end{enumerate}

The first Brezzi condition follows from the condition $\div \Sigma_h \subset V_h$. 
It is easy to see that the second follows from condition (2) above.
To construct elements which satisfy (1) and (2), we follow the constructive approach
of Arnold, Falk and Winther, \cite{Arnold2006,Arnold2007}, using 
discrete versions of the de Rham sequence. 
In addition to the spaces $\Sigma_h, V_h$ and $Q_h$, we also construct finite
element spaces $R_h \subset H(\div,\Omega,\mathbb{H})$ and 
$\Theta_h \subset H(\curl,\Omega,\mathbb{H})$ in such a way that the
following diagrams commute:
$$
\begin{CD}
H(\div,\Omega,\mathbb{H} ) 
@>\div >> L^2(\Omega, \mathbb{P} ) @>>> 0\\
                @VV{\Pi_{R_h}} V @VV{\Pi_{Q_h}}V  @.\\
R_h @>
\div>>
 Q_h @>>> 0,  
\end{CD}
$$
$$
\begin{CD}
H(\curl,\Omega,\mathbb{H}) @>\curl>> 
H(\div,\Omega,\mathbb{M}) 
@>\div >> L^2(\Omega,\mathbb{R}^n) @>>> 0\\
   @VV{\Pi_{\Theta_h}}V @VV{\Pi_{\Sigma_h}} V @VV{\Pi_{V_h}}V  @.\\
 \Theta_h  @>\curl>> \Sigma_h @>
\div>>
 V_h @>>> 0.  
\end{CD}
$$
We note that the commutativity of the far left side of the diagram above will not be used.
For a finite dimensional space $X_h$, $\Pi_{X_h}$ is a bounded projection operator. We recall that 
\begin{align}
\Pi_{X_h} v  = v, \ \forall \ v \in X_h. \label{proj2}
\end{align}
Next, we define an operator $S: C^{\infty}(\Omega,\mathbb{H} )\to C^{\infty}(\Omega,\mathbb{H})$ which 
connects the two diagrams above. 
In two dimension, $S$ is simply the identity operator, while in three dimension,
for $q=(q_{ij})_{i,j=1,\ldots,3}$, we define
\begin{equation} S(q) =
\begin{pmatrix}
q_{22} +q_{33} & -q_{21} & -q_{31}  \\
-q_{12} & q_{11} +q_{33} & -q_{32} \\
-q_{13} & -q_{23} & q_{11}+ q_{22} 
\end{pmatrix} .
\end{equation}
In that case, $S$ is also invertible with $S(q)=\tr(q)I -q^T, S^{-1}(q) = 1/2 \tr(q) I -q^T$, \cite{Falk08},
where $q^T$ denotes the transpose of $q$, $I$ is the $3 \times 3$
identity matrix and $\tr(q)$ denotes the trace of $q$.
The following fundamental relation holds in both dimension:
\begin{equation}
\as \curl q =-\div S(q). \label{fundamental}
\end{equation}
We summarize the elements of the
constructive approach of \cite{Arnold2006,Arnold2007} in the following theorem, the proof of which 
is reproduced below for convenience.
\begin{thm} \label{thmCond}
Under the commutativity assumptions  
\begin{align}
\Pi_{Q_h} \div q &= \div \Pi_{R_h} q, \, \forall q \in C^{\infty}(\Omega,\mathbb{H}),\label{com1}\\
\div \Pi_{\Sigma_h} \sigma &= \Pi_{V_h} \div \sigma, \,
\forall \, \sigma \in C^{\infty}(\Omega,\mathbb{M}), \label{com2}
\end{align}
and
\begin{align}
\Pi_{R_h}  S \Pi_{\Theta_h} S^{-1} &=\Pi_{R_h}, \label{ceq1} \\
||\Pi_{\Sigma_h}u||_{L^2} &\leq c ||u||_{H^1}, \ \forall \tau \in H^1(\Omega,\mathbb{M}), \label{ceq3} \\
||\curl \Pi_{\Theta_h}\rho||_{L^2} &\leq c ||\rho||_{H^1}, \ \forall \rho 
\in H^1(\Omega,\mathbb{H}). \label{ceq4} 
\end{align}
the second Brezzi condition holds.
\end{thm}
\begin{proof}
By elliptic regularity, given $v \in V_h,
\exists \, \eta \in H^1(\Omega,\mathbb{M})$ such that
\begin{equation}
\div \eta = v \quad \text{and} \quad ||\eta||_{H^1} \leq ||v||_{L^2}. \label{ceq5}
\end{equation} 
Given $q \in Q_h \subset L^2(\Omega,\mathbb{P})$, there exists 
$ \phi \in H^1(\Omega,\mathbb{H})$ such that
\begin{equation}
\div \phi = q - \Pi_{Q_h} \as \Pi_{\Sigma_h} \eta \ 
\text{and} \ ||\phi||_{H^1} \leq C ||q - \Pi_{Q_h} \as \Pi_{\Sigma_h} \eta ||_{L^2}.
\label{ceq6}
\end{equation} 
We set $\tau= \Pi_{\Sigma_h} \eta + \curl \Pi_{\Theta_h}S^{-1} \phi$ and by \eqref{com2}
and \eqref{proj2} we have
$$
\div \tau = \div \Pi_{\Sigma_h} \eta = \Pi_{V_h} \div \eta =\Pi_{V_h} v = v. 
$$
By \eqref{fundamental} and \eqref{com1} it follows that
$$
\Pi_{Q_h} \as \curl q = \Pi_{Q_h} \div S q = \div \Pi_{R_h} S q,
$$
We therefore have using \eqref{ceq1}, \eqref{com1} and \eqref{proj2},
\begin{align*}
\Pi_{Q_h} \as \tau &= \Pi_{Q_h} \as \Pi_{\Sigma_h} \eta + 
\Pi_{Q_h} \as \curl \Pi_{\Theta_h} S^{-1} \phi \\
& = \Pi_{Q_h} \as \Pi_{\Sigma_h} \eta + \div \Pi_{R_h} S \Pi_{\Theta_h} S^{-1} \phi \\
& = \Pi_{Q_h} \as \Pi_{\Sigma_h} \eta + \div \Pi_{R_h} \phi \\
& = \Pi_{Q_h} \as \Pi_{\Sigma_h} \eta + \Pi_{Q_h} \div \phi \\
& = \Pi_{Q_h} \as \Pi_{\Sigma_h} \eta + \Pi_{Q_h} q - \Pi_{Q_h} \as \Pi_{\Sigma_h} \eta\\
& = q.
\end{align*}
It remains to prove the inequality \eqref{divrel}. We have by \eqref{ceq5} and  \eqref{ceq3}
$$
||\Pi_{\Sigma_h} \eta||_{L^2} \leq C ||\eta||_{H^1} \leq C ||v||_{L^2},
$$ 
and by \eqref{ceq5}, \eqref{proj2}, \eqref{ceq5}, \eqref{ceq3} and \eqref{ceq6}
\begin{align*}
||\curl \Pi_{\Theta_h} S^{-1} \phi ||_{L^2} & \leq ||S^{-1}\phi||_{H^1} \leq C ||\phi||_{H^1} \leq 
||q - \Pi_{Q_h} \as \Pi_{\Sigma_h} \eta ||_{L^2} \\
& \leq C ( ||q||_{L^2} + ||\as \Pi_{\Sigma_h} \eta ||_{L^2}) 
\leq C( ||q||_{L^2} + ||\eta||_{H^1} )\\
& \leq C ( ||q||_{L^2} + ||v||_{L^2}).
\end{align*}
It follows that $ ||\tau||_{L^2} = ||\Pi_{\Sigma_h} \eta + \curl \Pi_{\Theta_h} \phi||_{L^2}
\leq C (||q||_{L^2}+||v||_{L^2})$. Since $\div \tau=v$, this proves the result.
\end{proof}

Let $\mathcal{T}_h$ denote a conforming partition of $\Omega$
 into rectangles of diameter bounded by $h$, which
is quasi-uniform in the sense that the aspect ratio of the rectangles is bounded by a fixed
constant. Let $\hat{R}=[0,1]^n$ be the reference rectangle and let $F: \hat{R} 
\to R$
be an affine mapping onto $R$, $F(\hat{x}) = B \hat{x} +b$, with $b \in \mathbb{R}^n$ 
and $B$ a $n \times n$ diagonal matrix.
Our goal in the next section is to construct spaces
$\Sigma_h, V_h$ and $\Theta_h$ such that the conditions
of Theorem \eqref{thmCond} hold. 
If $(\sigma,u,p)$ denotes the solution of  problem \eqref{elasSys} and
$(\sigma_h,u_h,p_h) \in \Sigma_h \times V_h \times \Theta_h$ is
the solution of the associated discrete system, the optimality condition
\begin{align}
\begin{split}\label{optimal}
||\sigma-\sigma_h||_{H(\div)} + ||u-u_h||_{L^2} & +||\gamma-\gamma_h||_{L^2}
 \leq C \, \text{inf}_{\tau_h \in \Sigma_h,
v_h \in V_h, \rho_h \in Q_h} \\& (||\sigma-\tau_h||_{H(\div)}  
+||u-v_h||_{L^2} +||\gamma-\rho_h||_{L^2}),
\end{split}
\end{align}
holds.

As with \cite{Arnold2006,Arnold2008,Falk08}, the following refined error estimates hold
\begin{align*}
||\sigma-\sigma_h||_{H(\div)} + ||u_h-\Pi_{V_h} u||_{L^2} & +||\gamma-\gamma_h||_{L^2}
 \leq C (||\sigma-\Pi_{\Sigma_h}\sigma||_{H(\div)}+||\gamma-\Pi_{Q_h}\gamma||_{L^2}), \\
||u-u_h||_{L^2}  \leq C (||\sigma &-\Pi_{\Sigma_h}\sigma||_{H(\div)} +||\gamma-\Pi_{Q_h}\gamma||_{L^2}
+||u -\Pi_{V_h} u||_{L^2}), \\ 
||\div(\sigma-\sigma_h)||_{L^2} & = ||\div \sigma - \Pi_{V_h} \div \sigma ||_{L^2}.
\end{align*}

\section{Two dimensional elements}
We recall the lowest order BDM element,
\begin{equation}
BDM_1(K) = \{ \, q \,| \, q=p_1(x,y)+r \curl(x^{2}y) +s \curl(x y^{2}),
p_1 \in \mathcal{P}_1 \times \mathcal{P}_1 \, \}, 
\end{equation}
and an element $q \in BDM_1(K)$ is uniquely determined by the conditions
$\int_e q \cdot n \, p_1 \ud s, \, \text{for each edge $e$ of $K$}, \, \forall
\, p_1 \in
\mathcal{P}_1(e) $.

We choose $V_h=\mathcal{P}_0(\mathcal{T}_h)$,
$Q_h=\mathcal{P}_0(\mathcal{T}_h)$,
with degrees of freedom the value at an interior point in each element $K$ and
$$
\Sigma_K = \{\, \tau, \tau(x,y) \in \mathbb{M}, (\tau_{i 1}, \tau_{i 2}) \in BDM_1(K), i=1,2 \,\}.
$$
A matrix field $\tau \in \Sigma_K$ is uniquely determined by the first two
moments of $\tau n$ on each edge, ($2\times 2\times 4=16$ degrees of freedom). The 
stress field space $\Sigma_h$ is therefore 
the space of matrix fields which belong piecewise to $\Sigma_K$ and have normal components 
which are continuous across mesh edges. 

We will also need the serendipity finite element space $S_h$, defined on a single element $K$ by
\begin{equation*}
S_K = \mathcal{P}_2(K) + \text{span}\{x^2 y, x y^2\},
\end{equation*}
and with degrees of freedom for $q \in S_K$
\begin{enumerate}
\item the values of $q$ at the vertices (4 degrees of freedom),
\item the average of $q$ on each edge (4 degrees of freedom).
\end{enumerate}
It is not difficult to check that the sequence
$$
\begin{CD}
0 @ >>> 
\mathbb{R} 
@>\subset>> 
S_K
@ > \curl >>
BDM_1(K)
@>\div >> 
\mathcal{P}_0(K)
@>>> 0.
\end{CD}
$$
is exact. One checks that each space is mapped in the one that follows. Then one notes that the alternating sum of the dimensions
is zero and that the polynomial de Rham sequence is exact.

We therefore define the space $\Theta_h$ as follows: on each element $K$,
$\Theta_K = S_K \times S_K$ and the space $\Theta_h$ 
is the space of vector fields which belong piecewise to $\Theta_K$ and 
are continuous across mesh edges.

Finally we take for $R_h$ the lowest order Raviart-Thomas element, i.e.
$R_h=RT_0(\mathcal{T}_h)$. We recall that
$ RT_0(K) = \mathcal{P}_{1,0}(K) \times \mathcal{P}_{0,1}(K)$ with degrees of freedom
the average of the normal component of $q \in RT_0(K)$ on each edge. 

The projection operator $ \Pi_{\Sigma_h} $ is taken as the canonical interpolation operator
and defined  by
\begin{align*}
\int_e  \Pi_{\Sigma_h}(\sigma) n\cdot q \ud s = \int_e \sigma n\cdot q \ud s, \quad
\text{ for all edges} \ e  \ \text{and for all} \ 
q \in \mathcal{P}_1(e) \times \mathcal{P}_1(e).
\end{align*}
Similarly we define $\Pi_{R_h}$  by
\begin{align*}
\int_e  \Pi_{R_h}(q) \cdot n \ud s = \int_e q \cdot n \ud s, \quad
\text{ for all edges } \ \ e. 
\end{align*}
It remains to define the interpolation operator $\Pi_{\Theta_h}$. For this we first define
$ \Pi^0_K: H^1(K,\mathbb{R}^2) \to \Theta_K$ by
\begin{align*}
\Pi^0_K \psi(v) & =0 \quad \text{for each vertex} \ v \ \text{of} \ K,\\
\int_e \Pi^0_K \psi(s) \ud s & = \int_e \psi(s) \ud s\quad \text{for each edge} \ e
\subset \partial K,
\end{align*}
and $ \Pi^0_h: H^1(\Omega,\mathbb{R}^2) \to \Theta_h$ by 
$(\Pi^0_h \tau)|_K = \Pi^0_K \tau$.
Next, let $L_h$ be a Clement interpolation operator \cite{Bernardi98,Cl'ement1975} which
maps $L^2(\Omega,\mathbb{R})$ into
$$
\{\, \theta_h \in C^0(\bar{\Omega})\, | \, \theta_{h|K} \in \mathcal{P}_{1,1},
\forall K \in \mathcal{T}_h \, \},
$$
and denote as well by $L_h$ the corresponding operator which maps 
$L^2(\Omega,\mathbb{R}^2)$
into the subspace $\Theta_h$ of continuous vector fields whose components 
are piecewise in $\mathcal{P}_{1,1}$.
We have
\begin{equation}
\|L_h \tau - \tau \|_j \leq c h^{m-j} \|\tau\|_m, \quad 0\leq j \leq 1,
\quad j \leq m 
\leq 2, \label{rhbound}
\end{equation}
with $c$ independent of $h$. We define our interpolation operator $\Pi_{\Theta_h}$ by
\begin{equation}
\Pi_{\Theta_h} = \Pi^0_h (I - L_h) + L_h. \label{opdef}
\end{equation}
\begin{thm}
For the triple $(\Sigma_h, V_h,\Theta_h)$ the conditions
of Theorem \eqref{thmCond} hold and we have the optimality condition \eqref{optimal}.
Moreover if $\sigma$ and $u$ are sufficiently smooth,
\begin{align}
||\sigma-\sigma_h||_{H(\div)} + ||u-u_h||_{L^2}  +||\gamma-\gamma_h||_{L^2}
 \leq C \, h ||u||_3. \label{lowest1}
\end{align}
\end{thm}
\begin{proof}
Let $q \in C^{\infty}(\Omega,\mathbb{R}^2)$. We have using the definition
of $\Pi_{R_h}$ and Green's theorem,
\begin{align*}
\int_{\Omega} \div \Pi_{R_h} q \ud x & = \sum_K \int_K \div \Pi_{R_h} q \ud x 
= \sum_K \int_{\partial K} \Pi_{R_h} q \cdot n \ud s\\
&  = \sum_K \int_{\partial K} q \cdot n \ud s = \int_{\Omega} \div q,
\end{align*}
which proves \eqref{com1}.

Next, let $\sigma \in C^{\infty}(\Omega,\mathbb{M})$. Again using the definition
of $\Pi_{\Sigma_h}$ and Green's theorem,
\begin{align*}
\int_{\Omega} \div \sigma - \div \Pi_{\Sigma_h} \sigma \ud x & =
\sum_K \int_{K} \div (\sigma - \Pi_{\Sigma_h} \sigma) \ud x
= \sum_K \int_{\partial K} (\sigma - \Pi_{\Sigma_h} \sigma) n \ud s
=0,
\end{align*}
which proves \eqref{com2}.

For $q \in C^{\infty}(\Omega,\mathbb{R}^2)$, put $u= \Pi_{h}^0 q $. We have
using the definition
of $\Pi_{h}^0$ 
\begin{align*}
\int_e (u-q)\cdot n \ud s= \int_e ( \Pi_{h}^0 q -q) \cdot n \ud s = 0.
\end{align*}
It follows that $\Pi_{R_h}(u-q)=0$ i.e. $\Pi_{R_h} \Pi_{h}^0 q = \Pi_{R_h} q$. Finally, 
$\Pi_{R_h}  \Pi_{\Theta_h}=\Pi_{R_h} \Pi_{h}^0  (I-L_h)+\Pi_{R_h} L_h 
= \Pi_{R_h}(I-L_h)+ \Pi_{R_h}L_h= \Pi_{R_h} $,
that is \eqref{ceq1} holds.

By the trace theorem, one shows that $ (\Pi_{\Sigma_h})|_{\hat{K}} $ is bounded on
$H^1(\hat{K},\mathbb{M})$. Moreover if we define
for a matrix field $\hat{M}$, $P_F(\hat{M})(x)= 1/\text{det}(B) \hat{M}(\hat{x}) B^T
, x=F(\hat{x})$,
then it is not difficult to verify that
$P_F((\Pi_{\Sigma_h})|_{\hat{K}} \hat{\sigma})
=(\Pi_{\Sigma_h})|_{K}P_F \hat{\sigma}$, hence \eqref{ceq3} follows from a standard scaling argument.

Let $ \hat{\rho} \in H^1(\hat{K},\mathbb{R}^2)$. 
We define its Piola transform by $ P_F \hat{\rho} =(P_F \hat{\rho}_1,P_F \hat{\rho}_2)$
where for a scalar function $\hat{u}$, $P_F \hat{u}=\hat{u}\circ F^{-1}$.

Since $\hat{\curl} \Pi^0_{\hat{K}} \hat{\rho}\in 
\Sigma_{\hat{K}}$, 
$$
|| \hat{ \curl} \Pi^0_{\hat{K}} \hat{\rho}||_{L^2(\hat{T})} \leq C \sum_{\hat{e} \subset \partial \hat{K}} 
\sum_{i=0}^1 \bigg| \int_{\hat{e}} \hat{\curl} \Pi^0_{\hat{K}} \hat{\rho}\cdot \hat{n} \hat{s}^i \ud \hat{s} \bigg|,
$$ 
where $\hat{e}$ is an edge of $\partial \hat{K}$. 
Next, $\curl q \cdot n = \partial q/\partial s$ and using the definition of $\Pi^0_{\hat{K}}$,
\begin{align*}
\int_{\hat{e}} \hat{\curl} \Pi^0_{\hat{K}} \hat{\rho}\cdot \hat{n} \ud \hat{s} &
= \int_{\hat{e}} \frac{\partial}{\partial \hat{s}} \Pi^0_{\hat{K}} \hat{\rho}\ud \hat{s} =0\\
\int_{\hat{e}} \hat{\curl} \Pi^0_{\hat{K}} \hat{\rho}\cdot \hat{n} \ \hat{s} \ud \hat{s} &
= \int_{\hat{e}} \frac{\partial}{\partial s} ( \Pi^0_{\hat{K}} \hat{\rho} ) \hat{s} \ud \hat{s}
= -\int_{\hat{e}} \Pi^0_{\hat{K}} \hat{\rho} \ud \hat{s} = -\int_{\hat{e}} \hat{\rho} \ud \hat{s}.
\end{align*}
By the trace theorem, it follows that
$$
|| \hat{\curl} \Pi^0_{\hat{K}} \hat{\rho}||_{L^2(\hat{T})} \leq C ||\hat{\rho}||_{1,\hat{T}},
$$
and scaling to an arbitrary rectangle $K$, we get
$$
||\curl \Pi^0_K \rho||_{L^2(K)} \leq C (h^{-1} |\rho|_{0,K} + C |\rho|_{1,K}).
$$
We therefore have
\begin{align*}
||\curl \Pi_{\Theta_h} \rho||_{L^2} & \leq 
||\curl \Pi^0_{h}(I-L_h)\rho||_{L^2} + ||\curl L_h \rho||_{L^2} \\
& \leq c (h^{-1} ||(I-L_h)\rho||_{L^2} + ||(I-L_h)\rho||_{H^1})  + c ||L_h \rho||_{H^1} \\
&  \leq c ||\rho||_{H^1}, 
\end{align*}
that is \eqref{ceq4} holds. Since $\div \Sigma_h \subset V_h$, the Brezzi 
conditions hold and the error estimates follow from
the optimality error estimate from the theory of
mixed methods, properties of the canonical interpolation operator for BDM elements,
\cite{Brezzi1991} p. 132, and error estimates of the $L^2$ projection operator.
\end{proof}

\subsection{Simplified element of low order}
Analogous to the simplified element of \cite{Arnold2006}, we can develop elements simpler than the lowest 
order BDM type elements. The key point is that for \eqref{ceq1} to hold, we only need $\Theta_h$
to have normal components continuous across edges. We start the construction by
taking as $\Theta_h$ the rectangular version of a space introduced by Fortin, \cite{Fortin81}
and \cite{Girault86} p. 153. The spaces $R_h$, $V_h$ and $Q_h$ are the same. To define
the space $\Theta_h$, let $i,j$ be the unit vectors in the $x$ and $y$ directions
respectively. We put
\begin{alignat*}{4}
p_1 & = -x(1-x)(1-y) \, i\\
p_2 & = -y(1-y)(1-x) \, j \\
p_3 & = x(1-x)y \, i\\
p_4 & = x y(1-y) \, j,
\end{alignat*} 
and define on each element $K$,
$$
\Theta_K = \mathcal{P}_{1,1}(K) \times \mathcal{P}_{1,1}(K) \oplus \text{span} \, \{ \,
p_1,p_2,p_3,p_4 \, \}
$$
with degrees of freedom
\begin{enumerate}
\item the values of $q$ at the vertices ($4\times 2=8$ degrees of freedom),
\item the average of $q \cdot n$ on each edge (4 degrees of freedom).
\end{enumerate}
The stress space $\tir{\Sigma}_K$ is defined as
$$
\begin{pmatrix}
\mathcal{P}_{1,0} (K) & \mathcal{P}_{0,1}(K)  \\
\mathcal{P}_{1,0}(K)  & \mathcal{P}_{0,1}(K)   
\end{pmatrix}
\oplus \text{span} \, \{ \,
\curl p_1, \curl p_2, \curl p_3, \curl p_4 \, \},
$$
where $\begin{pmatrix}
\mathcal{P}_{1,0} (K) & \mathcal{P}_{0,1}(K)  \\
\mathcal{P}_{1,0}(K)  & \mathcal{P}_{0,1}(K)   
\end{pmatrix}
$ is the space of matrix fields with components in the indicated spaces.
Explicitly, we have
$\curl p_1 = \begin{pmatrix}x(1-x)& (1-2x) (1-y)\\ 0 & 0 \end{pmatrix}$,
$\curl p_2 = \begin{pmatrix}0 & 0 \\(-1+2y)(1-x) & -y (1-y)  \end{pmatrix}$,
$\curl p_3 = \begin{pmatrix}x(1-x)& -(1-2x) y\\ 0 & 0 \end{pmatrix}$ and
$\curl p_4 = \begin{pmatrix}0 & 0 \\x (1-2y) & -y (1-y)  \end{pmatrix}$.

For $\tau \in \begin{pmatrix}
\mathcal{P}_{1,0} (K) & \mathcal{P}_{0,1}(K)  \\
\mathcal{P}_{1,0}(K)  & \mathcal{P}_{0,1}(K)   
\end{pmatrix}
$, $\tau n \in \mathcal{P}_{0}(e) \times \mathcal{P}_{0}(e) $ on each edge $e$ but 
$(\curl p_i) n \cdot t \in \mathcal{P}_{1}(e), i= 1, \ldots, 4$.
The following degrees of freedom are unisolvent:
\begin{enumerate}
\item $ \int_e \tau n \cdot n \ud s$ for each edge $e$
\item $  \int_e \tau n \cdot t \, p \ud s$ for each edge $e$ and $p \in \mathcal{P}_{1}(e)$.
\end{enumerate}
To see this, let $\tau = \eta + a_1 \curl p_1 + a_2 \curl p_2 + a_3 \curl p_3 + a_4 \curl p_4 \in \tir{\Sigma}_K$
such that all the above degrees of freedom vanish. Since the normal component of $(\tau_{i1},\tau_{i2}), i=1,2$ vanish
on each edge, we have
\begin{align*}
\tau_{i1} = x(1-x) c_{i1}, \tau_{i2} = y(1-y) c_{i2}, i=1,2, c_{i,j} \in \R, i,j=1,2.
\end{align*}
Since
\begin{align*}
\tau_{11} & = \eta_{11} + a_1 x(1-x) + a_3 x(1-x), \eta_{11} \in \P_{10}(K) \\
\tau_{12} & = \eta_{12} + a_1 (1-2x)(1-y) - a_3 (1-2x)y, \eta_{12} \in \P_{01}(K) \\
\tau_{21} & = \eta_{21} + a_2 (-1+2y)(1-x) - a_4 x(1-2y), \eta_{21} \in \P_{10}(K) \\
\tau_{22} & = \eta_{21} - a_4 y(1-y) - a_4 y(1-y), \eta_{22} \in \P_{01}(K),
\end{align*}
we conclude that $a_1=a_2=a_3=a_4=0$ and $\eta=0$, that is: $\tau=0$ and the claim follows.

From the approximation properties of the lowest order Raviart-Thomas element, the estimate
\eqref{lowest1} still holds.
\section{Three dimensional elements}
The de Rham complex in three dimensions is
$$
\begin{CD}
\mathbb{R} @ >\subset>> 
C^{\infty}(\Omega,\mathbb{R})
@>\grad>> C^{\infty}(\Omega,\mathbb{R}^3) 
@ > \curl >>
C^{\infty}(\Omega,\mathbb{R}^3)
@>\div >> 
C^{\infty}(\Omega,\mathbb{R})
@>>> 0.
\end{CD}
$$
We choose the following form of BDM elememt, \cite{Brezzi1991}, p.124
\begin{align*}
BDM_1(K) =  \mathcal{P}_{1}(K,\R^3)  +\curl \text{span} \bigg\{\, 
\begin{pmatrix}
0 \\0 \\ xy^2
\end{pmatrix},
\begin{pmatrix}
0\\0\\x^2y
\end{pmatrix},
\begin{pmatrix}
y^2z \\0 \\0
\end{pmatrix},
\begin{pmatrix}
y z^2 \\0 \\0
\end{pmatrix},
\begin{pmatrix}
0 \\ xz^2 \\0
\end{pmatrix},
\begin{pmatrix}
0 \\x^2z \\0
\end{pmatrix}
\, \bigg\}.
\end{align*}
Clearly $\div BDM_1(K)= \mathcal{P}_0(K)$. We define $V_K=\mathcal{P}_0(K)^3$ and
$$
\Sigma_K = \{\, \tau, \tau(x,y,z) \in \mathbb{M}, (\tau_{i 1}, \tau_{i 2}, \tau_{i 3}) \in BDM_1(K), i=1,2,3 \,\}.
$$
The degrees of freedom on $V_K$ are the values of each component at an interior point while a matrix
field $\tau$ in $\Sigma_K$ is uniquely determined by the moments of order 0 and 1 of $\tau n$ on each 
face ($3 \times 3 \times 6$ degrees of freedom).

We now define two spaces
$S_K$ and $U_K$ such that the sequence below is exact.
$$
\begin{CD}
\mathbb{R} @ >\subset>> 
S_K
@>\grad>> U_K
@ > \curl >>
BDM_1(K)
@>\div >> 
\P_0(K, \mathbb{R} )
@>>> 0
\end{CD}.
$$
The space $S_K$ is not directly used in the construction but helped discover $U_K$.
We take the space $S_K$ as the three dimensional serendipity space of order 2 defined as
\begin{align*}
\begin{split}
S_K & =  \mathcal{P}_{2}(K,\R) + \text{span}   \{\, 
x^2y, x^2z, xy^2 ,xz^2, y^2z, yz^2, xyz,
x^2yz,xy^2z,xyz^2
\,\}, 
\end{split}
\end{align*}
with degrees of freedom
\begin{enumerate}
\item the values of $q \in S_K$ at the vertices (8 degrees of freedom),
\item the average of $q \in S_K$ on each edge (12 degrees of freedom).
\end{enumerate}
The unisolvency of these degrees of freedom is proven for example in \cite{Arnold2010d}.
We define the space $U_K$ as
\begin{align*}
\begin{split}
U_K & = 
\mathcal{P}_{1,1,1}(K,\mathbb{R}^3) +  \text{span}  \{\, y^2z,yz^2,y^2,z^2 \,\} \times  \text{span}  \{\,x^2z,xz^2,x^2,z^2 \,\}  \times  \text{span}  
\{\, x^2y,xy^2,x^2,y^2 \,\}, 
\end{split}
\end{align*}
with degrees of freedom for $u \in U_K$,
\begin{enumerate}
\item the first two moments of $u \cdot t$ on each edge, where
$t$ is a tangential vector to the edge ($12\times 2=24$ degrees of freedom),
\item the average of $u \wedge n$ on each face 
with unit outward normal $n$ ($6 \times 2 =12$ degrees of freedom).
\end{enumerate}

It is not very difficult to verify that
the sequence above is exact.
One checks that each space is mapped in the one that follows. Then one notes that the alternating sum of the dimensions
is zero and that the polynomial de Rham sequence is exact. We then only need to verify either that the kernel of the $\curl$
operator is the image of the $\grad$ operator or that the kernel of the $\div$ operator is the image of the $\curl$ operator.
We verify the last one. Let $u \in BDM_1(K)$ such that $\div u=0$. We write $u=w + \curl z, w \in \mathcal{P}_{1}(K,\R^3)$
and $z$ in the span of the extra monomials in the definition of $BDM_1(K)$. Note that $z \in U_K$ and $\div u = \div w=0$.
By the exactness of the polynomial de Rham sequence, $w=\curl a, a \in \mathcal{P}_{2}(K,\R^3)$. Since for
$\alpha, \beta, \gamma \in \R, \curl(\alpha x^2, \beta y^2, \gamma z^2)=0$, we may assume that $a \in U_K$ which completes the proof of the claim.

We can now describe the space $\Theta_h$ as
$$
\Theta_h = \{\, q, q(x,y,z) \in \mathbb{M}, (q_{i1},q_{i2},q_{i3}) \in U_h, i=1,2,3\,\}, 
$$
with the degrees of freedom for $q \in \Theta_h$
\begin{enumerate}
\item $\int_e q \, t \,s^i, i=0,1$ for each edge $e$, where
$t$ is a tangential vector to the edge ($12\times 2 \times 3=72$ degrees of freedom),
\item $\int_f q \wedge n \ud x_f$ for each face $f$ with unit outward normal $n$ 
($6 \times 2 \times 3=36$ degrees of freedom). For a matrix field $q$ with row vectors
$q_i, i=1,2,3$, $q \wedge n$ is defined
as the matrix field with rows $q_i \wedge n, i=1,2,3$.
\end{enumerate}

Next we define the space $Q_h$. We take $Q_K=\mathcal{P}_0(K)^3$ with degrees of freedom the values of each component at an interior point.

Finally we describe the space $R_h$ as 
$$
\{\, q, q(x,y,z) \in \mathbb{M}, (q_{i1},q_{i2},q_{i3})_{|_K} \in RT_0(K), i=1,2,3\,\}, 
$$
where
\begin{align*}
RT_0(K) = \mathcal{P}_{1,0,0}(K) \times \mathcal{P}_{0,1,0}(K) \times \mathcal{P}_{0,0,1}(K),
\end{align*}
is the lowest order Raviart-Thomas element in three dimensions with degrees of freedom the average of the normal component on each face, ($1 \times 1 \times 6$=6 degrees of freedom).

\subsubsection{Unisolvency}
The unisolvency of the degrees of freedom for $V_K$, $\Sigma_K$ and $S_K$ are well known. Similarly
unisolvency for the degrees of freedom of $R_h$ is immediate. We only study the case of $U_K$.
Let $v=(v_1,v_2,v_3) \in U_K$ and assume that all degrees of freedom vanish. We show that $v_1=0$.
On each edge $e$, $v \cdot t \in \P_1(e)$ and hence we get $v \cdot t=0$ on each edge.
This implies that on the face $z=0$ for example,
\begin{align*}
v_1 & = y(1-y) w_1, w_1 \in \P_{1,0}\\
v_2 & = x(1-x) w_2, w_2 \in \P_{0,1}.
\end{align*}
However, if $w_1$ has a linear term in $x$, $xy^2$ would be the highest degree monomial in $v_1$. We conclude that
$w_1$ is constant. The face degrees of freedom imply that the average of $w_1$ vanish on the face $z=0$, that is: $w_1=0$.
Similarly $w_2=0$. We conclude that $v$ has expression
\begin{align*}
v_1 & = y(1-y)z (1-z) r_1, \\
v_2 & = x(1-x) z(1-z) r_2, \\
v_3 & = x(1-x) y(1-y) r_3,
\end{align*}
for constants $r_1, r_2$ and $r_3$ which must vanish given the form of the highest degree monomial in the expression of 
$v_i, i=1,2,3$.
\subsubsection{Definition of interpolation operators }
For $q \in C^{\infty}(\Omega,\mathbb{M})$, we define $\Pi_{R_h}$ by
$$
\int_{f} (\Pi_{R_h} q) n \ud x = \int_{f} q n \ud x, \quad \text{for all faces} \ f. 
$$
The interpolation operator $ \Pi_{\Sigma_h} $ 
is defined  by
\begin{align*}
\int_f  \Pi_{\Sigma_h}(\sigma) n\cdot q \ud s = \int_f \sigma n\cdot q \ud s, \quad
\text{ for all faces} \ f  \ \text{and for all} \ 
q \in \mathcal{P}_1(f) \times \mathcal{P}_1(f) \times \mathcal{P}_1(f).
\end{align*}
It remains to define the interpolation operator $\Pi_{\Theta_h}$. For this we first define
$ \Pi^0_K: H^1(K,\mathbb{M}) \to \Theta_K$ by
\begin{align*}
\int_e (\Pi^0_K q) \, t \, s^i \ud s & = 0, \, i=0,1 \quad \text{for each edge} \ e
\subset \partial K, \\
\int_f (\Pi^0_K q) \wedge n \ud x_f
& = \int_f q \wedge n \ud x_f, \quad \text{for each face} \ f
\subset \partial K
\end{align*}
and $ \Pi^0_h: H^1(\Omega,\mathbb{M}) \to \Theta_h$ by 
$(\Pi^0_h \tau)|_K = \Pi^0_K \tau$.
Next, let $L_h$ be a Clement interpolation operator \cite{Bernardi98,Cl'ement1975} which
maps $L^2(\Omega,\mathbb{R})$ into
$$
\{\, \theta_h \in C^0(\bar{\Omega})\, | \, \theta_{h|K} \in \mathcal{P}_{1,1,1},
\forall K \in \mathcal{T}_h \, \},
$$
and denote as well by $L_h$ the corresponding operator which maps 
$L^2(\Omega,\mathbb{M})$
into the subspace of $\Theta_h$ of continuous matrix fields whose components 
are piecewise in $\mathcal{P}_{1,1,1}$.
We have
\begin{equation}
\|L_h \tau - \tau \|_j \leq c h^{m-j} \|\tau\|_m, \quad 0\leq j \leq 1,
\quad j \leq m 
\leq 2, \label{rhbound}
\end{equation}
with $c$ independent of $h$. We define our interpolation operator $\Pi_{\Theta_h}$ by
\begin{equation}
\Pi_{\Theta_h} = \Pi^0_h (I - L_h) + L_h. \label{opdef}
\end{equation}

\subsubsection{Commutativity and surjectivity assumptions}
The commutativity assumption \eqref{com1} and \eqref{com2} are proven
as in the 2D case. We verify the surjectivity assumption 
$\Pi_{R_h}  S \Pi_{\Theta_h}=\Pi_{R_h} S $. We first show that
$\Pi_{R_h}  S \Pi_{\Theta_h}=\Pi_{R_h} S $. For this
let $q \in C^{\infty}(\Omega, \mathbb{M})$, put $\omega = q -\Pi^0_h q$.
We need to show that $\Pi_{R_h} S \omega =0$, that is
$$
\int_f (S \omega)(x) n \ud x_f = 0, \quad \text{for each face} \ f.
$$
Since $\Pi^0_h w =0$,
$$
\int_f \omega \wedge n =0, \quad \text{for each face} \ f.
$$
Next for $q=(q_{ij})_{i,j=1,2,3}$,
$$
q \wedge n =
\begin{pmatrix}
 q_{13} n_1 -q_{11} n_3 & q_{11} n_2 - q_{12} n_1 & q_{12} n_3 - q_{13} n_2 \\
q_{23} n_1 -q_{21} n_3 & q_{21} n_2 - q_{22} n_1 & q_{22} n_3 - q_{23} n_2 \\
q_{33} n_1 -q_{31} n_3 & q_{31} n_2 - q_{32} n_1 & q_{32} n_3 - q_{33} n_2 
\end{pmatrix},
$$
and
$$
(S q) n = \begin{pmatrix}
q_{22} n_1+ q_{33} n_1- q_{21} n_2- q_{31} n_3 \\
-q_{12} n_1+  q_{11} n_2+ q_{33} n_2- q_{32} n_3 \\
-q_{13} n_1- q_{23} n_2+ q_{11} n_3+ q_{22} n_3 
\end{pmatrix} = 
\begin{pmatrix}
-(q \wedge n)_{22}+ (q \wedge n)_{31} \\
(q \wedge n)_{12} - (q \wedge n)_{33} \\
-(q \wedge n)_{11} + (q \wedge n)_{23}
\end{pmatrix}.
$$
This shows that $\int_f \omega \wedge n =0$ implies $\int_f (S \omega) n =0$ and the result
follows using the definition of $\Pi_h$.

We notice that for $q \in \Theta_h$, for the surjectivity assumption to hold,
the following degrees of freedom were not used: $ \int_f q_{12} n_3 - q_{13} n_2  \ud x_f
= \int_f (q \wedge n)_{13}, 
\int_f q_{23} n_1 -q_{21} n_3 \ud x_f=\int_f (q \wedge n)_{12} , 
\int_f q_{31} n_2 - q_{32} n_1  \ud x_f=\int_f (q \wedge n)_{32}$.
However since the faces of a rectangle are parallel to the axes, one of these degrees of 
freedom is 
identically zero for each face, hence two degrees of freedom per face are unnecessary.
\subsubsection{Boundedness of the interpolation operators}
By the trace theorem, one shows that $ (\Pi_{\Sigma_h})|_{\hat{K}} $ is bounded on
$H^1(\hat{K},\mathbb{M})$. Moreover if we define
for a matrix field $\hat{M}$, $P_F(\hat{M})(x)= 1/\text{det}(B) \hat{M}(\hat{x}) B^T
, x=F(\hat{x})$,
then it is not difficult to verify that
$P_F((\Pi_{\Sigma_h})|_{\hat{K}} \hat{\sigma})
=(\Pi_{\Sigma_h})|_{K}P_F \hat{\sigma}$, hence \eqref{ceq3} follows from a standard scaling argument.

Let $ \hat{\rho} \in H^1(\hat{K},\mathbb{R}^3)$. 
We define its Piola transform by $ P_F \hat{\rho} =(P_F \hat{\rho}_1, P_F \hat{\rho}_2,P_F \hat{\rho}_3)  $
where for a scalar function $\hat{u}$, $P_F \hat{u}=\hat{u}\circ F^{-1}$.

Since $\hat{\curl} \Pi^0_{\hat{K}} \hat{\rho}\in 
\Sigma_{\hat{K}}$, 
$$
|| \hat{ \curl} \Pi^0_{\hat{K}} \hat{\rho}||_{L^2(\hat{T})} \leq C \sum_{\hat{f} \subset \partial \hat{K}} 
\sum_{i=0}^1 \bigg| \int_{\hat{f}} \hat{\curl} \Pi^0_{\hat{K}} \hat{\rho}\cdot \hat{n} \hat{s}^i \ud \hat{s} \bigg|,
$$ 
where $\hat{f}$ is a face of $\partial \hat{K}$.
Next, using the definition of $\Pi^0_{\hat{K}}$, for $q \in \mathcal{P}_{1,1}(f)
\times \mathcal{P}_{1,1}(f) \times \mathcal{P}_{1,1}(f)$,
\begin{align*}
\int_{\hat{f}} (\hat{\curl}(\Pi^0_{\hat{K}} \hat{\rho} ) \hat{n})\cdot q \ud x_f = 
\int_{\hat{f}} (\Pi^0_{\hat{K}} \hat{\rho}) \wedge \hat{n} \nabla q \ud x_f =
\int_{\hat{f}} \hat{\rho}  \wedge \hat{n} \nabla q \ud x_f .
\end{align*}
By the trace theorem, it follows that
$$
|| \hat{\curl} \Pi^0_{\hat{K}} \hat{\rho}||_{L^2(\hat{T})} \leq C || \hat{\rho} ||_{1,\hat{T}},
$$
and scaling to an arbitrary rectangle $K$, we get
$$
||\curl \Pi^0_K \rho||_{L^2(K)} \leq C (h^{-1} |\rho|_{0,K} + C |\rho|_{1,K}).
$$
We therefore have
\begin{align*}
||\curl \Pi_{\Theta_h} \rho||_{L^2} & \leq 
||\curl \Pi^0_{h}(I-L_h)\rho||_{L^2} + ||\curl L_h \rho||_{L^2} \\
& \leq c (h^{-1} ||(I-L_h)\rho||_{L^2} + ||(I-L_h)\rho||_{H^1})  + c ||L_h \rho||_{H^1} \\
&  \leq c ||\rho||_{H^1}, 
\end{align*}
that is \eqref{ceq4} holds. Since $\div \Sigma_h \subset V_h$, the Brezzi 
conditions hold. From
the optimality error estimate from the theory of
mixed methods \eqref{optimal}, properties of the canonical interpolation operator for BDM elements,
\cite{Brezzi1991} p. 132, and error estimates of the $L^2$ projection operator,
we have the following error estimate.
\begin{thm}
For the triple $(\Sigma_h, V_h,\Theta_h)$ the conditions
of Theorem \eqref{thmCond} hold and we have the optimality condition \eqref{optimal}.
Moreover if $\sigma$ and $u$ are sufficiently smooth,
\begin{align}
||\sigma-\sigma_h||_{H(\div)} + ||u-u_h||_{L^2}  +||\gamma-\gamma_h||_{L^2}
 \leq C \, h ||u||_3. \label{lowest}
\end{align}
\end{thm}
\section{Higher order elements}
Except the simplified element in two dimension, the elements we have described do not have optimal rate of convergence
for the stress. It does not seem possible to simplify the three dimensional element using the framework described here.
In two dimension, for higher order approximation, $H(\div)$ elements can be constructed based on the sequence,
$$
\begin{CD}
0 @ >>> 
\mathbb{R} 
@>\subset>> 
\mathcal{P}_{k+1,k+1}
@ > \curl >>
\mathcal{P}_{k+1,k} \times \mathcal{P}_{k,k+1}
@>\div >> 
\mathcal{P}_{k,k}
@>>> 0.
\end{CD}
$$
Take $V_h$ to be the space of piecewise continuous 
vector fields
which belong locally to $P_{k,k}(K) \times P_{k,k}(K)$, $Q_h$
the space of piecewise continuous functions
which belong locally to
$Q_K=P_{k-1,k-1}(K)$ and 
$\Sigma_K = \{\, \tau \in \mathbb{M}, (\tau_{i 1}, \tau_{i 2}) \in 
\mathcal{P}_{k+1,k} \times \mathcal{P}_{k,k+1}, i=1,2 \,\}$
with degrees of freedom
\begin{enumerate}
\item $\int_e \tau n \cdot p_k \ud s, \qquad \text{for each edge $e$ of $K$}, \ \forall
\, p_k \in
\mathcal{P}_k(e) $,
\item $\int_K \tau : \phi \ud x, \qquad \forall \, \phi \in 
\begin{pmatrix}
\mathcal{P}_{k,k-1} (K) & \mathcal{P}_{k-1,k}(K)  \\
\mathcal{P}_{k,k-1}(K)  & \mathcal{P}_{k-1,k}(K)   
\end{pmatrix}$,
\end{enumerate} 
for $k \geq 1$. The space $R_h$ is taken to be the Raviart-Thomas space of order $k-1$ and finally the space 
$\Theta_h$ is the space of continuous vector fields with components in 
$\mathcal{P}_{k+1,k+1}(K)$ on each element $K$.
Again, there one does not have optimal convergence rate for the stress.
We leave the details of the three dimensional analogue to 
the interested reader.
\bibliographystyle{amsplain}
\bibliography{SplineMonge03}

\providecommand{\bysame}{\leavevmode\hbox to3em{\hrulefill}\thinspace}
\providecommand{\MR}{\relax\ifhmode\unskip\space\fi MR }
\providecommand{\MRhref}[2]{%
  \href{http://www.ams.org/mathscinet-getitem?mr=#1}{#2}
}
\providecommand{\href}[2]{#2}
\begin{thebibliography}{10}

\bibitem{Adams2005}
Scot Adams and Bernardo Cockburn, \emph{A mixed finite element method for
  elasticity in three dimensions}, J. Sci. Comput. \textbf{25} (2005), no.~3,
  515--521.

\bibitem{Amara1979}
M.~Amara and J.~M. Thomas, \emph{Equilibrium finite elements for the linear
  elastic problem}, Numer. Math. \textbf{33} (1979), no.~4, 367--383.

\bibitem{Arnold2005}
Douglas~N. Arnold and Gerard Awanou, \emph{Rectangular mixed finite elements
  for elasticity}, Math. Models Methods Appl. Sci. \textbf{15} (2005), no.~9,
  1417--1429.

\bibitem{Arnold2010d}
\bysame, \emph{The serendipity family of finite elements}, To appear in
  Foundations of Computational Mathematics, November 2010., 2010.

\bibitem{Arnold2008}
Douglas~N. Arnold, Gerard Awanou, and Ragnar Winther, \emph{Finite elements for
  symmetric tensors in three dimensions}, Math. Comp. \textbf{77} (2008),
  no.~263, 1229--1251.

\bibitem{Arnold1984a}
Douglas~N. Arnold, Franco Brezzi, and {Jr} {Douglas, J.}, \emph{P{EERS}: a new
  mixed finite element for plane elasticity}, Japan J. Appl. Math. \textbf{1}
  (1984), no.~2, 347--367.

\bibitem{Arnold2006}
Douglas~N. Arnold, Richard~S. Falk, and Ragnar Winther, \emph{Finite element
  exterior calculus, homological techniques, and applications}, Acta Numer.
  \textbf{15} (2006), 1--155.

\bibitem{Arnold2007b}
\bysame, \emph{Finite element differential forms}, Proc. Appl. Math. Mech.
  \textbf{7} (2007), 1021901--1021902.

\bibitem{Arnold2007}
\bysame, \emph{Mixed finite element methods for linear elasticity with weakly
  imposed symmetry}, Math. Comp. \textbf{76} (2007), no.~260, 1699--1723
  (electronic).

\bibitem{Arnold2002}
Douglas~N. Arnold and Ragnar Winther, \emph{Mixed finite elements for
  elasticity}, Numer. Math. \textbf{92} (2002), no.~3, 401--419.

\bibitem{Arnold2003}
\bysame, \emph{Nonconforming mixed elements for elasticity}, Math. Models
  Methods Appl. Sci. \textbf{13} (2003), no.~3, 295--307, Dedicated to Jim
  Douglas, Jr. on the occasion of his 75th birthday.

\bibitem{Awanou10b}
Gerard Awanou, \emph{Symmetric matrix fields in the finite element method},
  Symmetry \textbf{2} (2010), 1375--1389.

\bibitem{Awanou10}
Gerard Awanou, \emph{Two remarks on rectangular mixed finite elements for
  elasticity in three dimensions}, Submitted, 2010.

\bibitem{Bernardi98}
C.~Bernardi and V.~Girault, \emph{A local regularization operator for
  triangular and quadrilateral finite elements}, SIAM J. Numer. Anal.
  \textbf{35} (1998), no.~5, 1893--1916 (electronic).

\bibitem{Falk08}
Daniele Boffi, Franco Brezzi, Leszek~F. Demkowicz, Ricardo~G. Dur{\'a}n,
  Richard~S. Falk, and Michel Fortin, \emph{Mixed finite elements,
  compatibility conditions, and applications}, Lecture Notes in Mathematics,
  vol. 1939, Springer-Verlag, Berlin, 2008, Lectures given at the C.I.M.E.
  Summer School held in Cetraro, June 26--July 1, 2006, Edited by Boffi and
  Lucia Gastaldi.

\bibitem{Brezzi1991}
Franco Brezzi and Michel Fortin, \emph{Mixed and hybrid finite element
  methods}, Springer Series in Computational Mathematics, vol.~15,
  Springer-Verlag, New York, 1991.

\bibitem{Chen2010}
Shao-Chun Chen and Ya-Na Yang, \emph{Conforming rectangular mixed finite
  elements for elasticity}, J. of Scientific Computing (2010), Submitted.

\bibitem{Cl'ement1975}
Ph. Cl{\'e}ment, \emph{Approximation by finite element functions using local
  regularization}, Rev. Fran\c caise Automat. Informat. Recherche
  Op\'erationnelle S\'er., RAIRO Analyse Num\'erique \textbf{9} (1975),
  no.~R-2, 77--84.

\bibitem{Guzman10d}
Bernardo Cockburn, Jayadeep Gopalakrishnan, and Johnny Guzm{\'a}n, \emph{A new
  elasticity element made for enforcing weak stress symmetry}, Math. Comp.
  \textbf{79} (2010), no.~271, 1331--1349.

\bibitem{Fortin81}
Michel Fortin, \emph{Old and new finite elements for incompressible flows},
  Internat. J. Numer. Methods Fluids \textbf{1} (1981), no.~4, 347--364.

\bibitem{Girault86}
Vivette Girault and Pierre-Arnaud Raviart, \emph{Finite element methods for
  {N}avier-{S}tokes equations}, Springer Series in Computational Mathematics,
  vol.~5, Springer-Verlag, Berlin, 1986, Theory and algorithms.

\bibitem{Guzman10b}
J.~Gopalakrishnan and J.~Guzm\'an, \emph{A second elasticity element using the
  matrix bubble}, Submitted, 2010.

\bibitem{Guzman10c}
J.~Guzm{\'a}n, \emph{A unified analysis of several mixed methods for elasticity
  with weak stress symmetry}, J. Sci. Comput. \textbf{44} (2010), no.~2,
  156--169.

\bibitem{Morley1989}
Mary~E. Morley, \emph{A family of mixed finite elements for linear elasticity},
  Numer. Math. \textbf{55} (1989), no.~6, 633--666.

\bibitem{Nicaise2008}
Serge Nicaise, Katharina Witowski, and Barbara~I. Wohlmuth, \emph{An a
  posteriori error estimator for the {L}am\'e equation based on equilibrated
  fluxes}, IMA J. Numer. Anal. \textbf{28} (2008), no.~2, 331--353.

\bibitem{Stenberg1986}
R.~Stenberg, \emph{On the construction of optimal mixed finite element methods
  for the linear elasticity problem}, Numer. Math. \textbf{48} (1986), no.~4,
  447--462.

\bibitem{Stenberg1988a}
\bysame, \emph{A family of mixed finite elements for the elasticity problem},
  Numer. Math. \textbf{53} (1988), no.~5, 513--538.

\bibitem{Stenberg1988}
\bysame, \emph{Two low-order mixed methods for the elasticity problem}, The
  mathematics of finite elements and applications, VI (Uxbridge, 1987),
  Academic Press, London, 1988, pp.~271--280.

\end{thebibliography}

\end{document}